\documentclass{arti}

\usepackage{graphics} 
\usepackage{epsfig} 
\usepackage{caption}
\usepackage{float}

\title{A Note on the Discrete Jordan Curve Theorem (Revised)}
\date{}
\author                        {  Li Chen\\
			         Department of Computer Science and Information Technology\\
				University of the District of Columbia	  \\
								 Washington, DC 20008\\
								 moorechen@yahoo.com
}

\begin{document}
\maketitle
\baselineskip=12pt

\centerline{\large {Abstract}}

According to a general definition of
discrete curves, surfaces, and manifolds ({\it Li Chen,
``Generalized discrete object tracking algorithms and
implementations,'' In Melter, Wu, and Latecki ed, Vision Geometry VI,
SPIE Vol. 3168, pp 184-195, 1997. }). This paper focuses on the Jordan curve
theorem in 2D discrete spaces. The Jordan curve theorem says that
a (simply) closed curve separates a simply connected surface into
two components.  Based on the definition of discrete surfaces,
we give three reasonable definitions of simply connected spaces.
Theoretically, these three definition shall be equivalent. We
have proved the Jordan curve theorem under the third definition
of simply connected spaces. The Jordan theorem shows the relationship
among an object, its boundary, and its outside area.

After the publication of the first version of the paper ({\it  L. Chen,  Note on the discrete Jordan Curve Theorem.
In: SPIE Conf. on Vision Geometry VIII, vol. 3811, pp. 82-94. SPIE, Bellingham (1999).} ), we found some
statements in the original proof of the Jordan Curve Theorem were not explained well.
One case was not proven in details. In this revision, we added
two more minor definitions and make the proof more solid and sound when it is needed for embedding
a discrete surface into a Euclidean space.

In this revision, we also proved  that the third definition
of simply connected spaces equivalent to the second definition
of simply connected spaces.


\section {Introduction}

In 1999,
The author declared that we have proved the discrete Jordan Curve Theorem for 2D discrete manifolds,
which should be a generalization of the proof of
the theorem in planar graphs made by W.T. Tutte in 1979$^{23}$.  However, when the author present the same result in the
book ``Discrete Surfaces and Manifolds'',
SP Computing,  in 2004$^{25}$. The author seems found that the proof in 1999 was not very  sound comparing to Tutte's proof.
Recently, researchers still show
considerable interests in Jordan curve theorem due to the fact of trustable mathematical proof and formalized proof using computers$^{24}$.
The author get back to review the original paper  by O. Veblen in 1905$^{23}$, this paper was regraded as the first correct proof of this fundamental theorem.
It was so interesting that Veblen's proof was not solidly sound neither if using a strict standard. This is because that Veblen made many definitions for some concepts
instead of proving them. It was very obvious that Veblen's definitions are very reasonable.

(It was shown that, for a very basic mathematics theorem, the proof may need to involve some reasonable axiom-like definitions in order to make a proof. )

In this revision, we will add two reasonable definitions, and they are much less than what Veblen did. We then will use the same proving strategy in our 1999 paper to
give a solid proof of  the discrete Jordan Curve Theorem. In fact,  the discrete Jordan Curve Theorem can be applied to prove the continuous version of Jordan
curve theorem in plane in a natural way. Our proof still maintain its pure discrete forms.

A simple closed curve will separate a plane into two disconnected components. This obvious fact needs a proof. C. Jordan gave a first proof. But many
mathematicians do not accept his proof. In 1905, Vablen gave a proof of this theorem. Vablen's proof not only provided a solid proof but also set
a standard for mathematics for what a proof should be.

In Veblen's proof, he first defined what a curve is? A curve must be defined. We know today, a 2D curve is a continuous mapping from [0,1] to a 2D Euclidean space.
This is even harder to define since we need to first define a function or a mapping. Draw a curve in the plane, can you get a mapping easily? No! 1-1 mapping for
[0,1] to the plane is not constructible. Veblen rather define the curve directly.

This is one of the motivation we want to use discrete curves for computers. Computers always want an object to be programable or even finite . that is possible to be put in
to the computer memory or disk space.

In order to make this current revision focuses only on  the discrete Jordan Curve Theorem,  we will omit some materials from the original paper.

Most of the related research work deals with curves and surfaces in $\Sigma_{m}$,
a space containing all integer grid points.   We usually call such a
space a digital space. A general definition of discrete
objects such as curves, surfaces, and $n$-Manifold is considered based on
a general graph $G=(V,E)$ in $^{1,4}$. This note focuses on the Jordan curve
theorem in the general 2D discrete space.\\\\

\section { Discrete Surfaces and Discrete Manifolds}

A graph, $G$, consists of two sets $V$ and $E$. $V$ is a set of vertices,
and $E$ is a set of pairs of vertices, called edges. An edge is said to be incident
with the vertices it joins.  We assume $G=(V,E)$ is an undirected graph
in this paper.
We also assume  $G=(V,E)$ is a simple graph, meaning every pair of vertices
has only one edge which is incident to these two vertices,
and there is no $(a,a) \in E$ for any  $a\in V$.

If $(p,q)$ is in $E$, $p$ is said to be adjacent to $q$.
Let $p_0,p_1,...,p_{n-1},p_{n}$ be $n+1$ vertices in $V$.
If $(p_{i-1},p_{i})$ is in $E$ for all $i=1,...,n$,
then $\{p_0,p_1,...,p_{n-1},p_{n} \}$ is called a path.
If $p_0,p_1,...,p_{n-1},p_{n}$ are distinct vertices,
the path is called a simple path.
A simple path $\{p_0,p_1,...,p_{n-1},p_{n} \}$ is closed if
$(p_0, p_{n})$ is an edge in $E$. A closed path is also called a cycle.
Two vertices $p$ and $q$ are connected if there is a path
$\{p_0,p_1,...,p_{n-1},p_{n} \}$ such that $p_0=p$ and $p_n=q$.
$G$ is called connected if every pair of vertices in $G$ is connected.
We always assume $G$ is connected in this paper.

Let $S$ be a set. Assume $S'$ is a subset of $S$, denoted by $S'\subseteq S$.
If $S$ is not a subset of $S'$, then $S'$ is called a proper-subset
of $S$, denoted by $S'\subset S$ .
Suppose  $G'=(V',E')$ is a graph where $V'\subset V$ and $E'\subset E$ for graph
$G=(V,E)$.  We say $G'=(V',E')$  is a partial-graph of $G$.
If $E'$ consists of all edges in $G$ whose joining vertices are in $V'$,
then the partial-graph $G'=(V',E')$ is called a subgraph of $G$, denoted
by $G'\preceq G$. If $V'$ is a proper-subset of $V$, then we denote $G'\prec G$.
We note that for a certain subset $V'$ of $V$,
the subgraph $G'$ with vertices $V'$ is uniquely defined.\\

{\bf 2.1 Basic Cells: Point-cells, Line-cells, and Surface-cells}

{\bf Definition 2.1} Each element of $V$ is  called a point-cell, 0-cell, or point.
Each element of $E$ is called a line-cell, or 1-cell.

Because $G$ is simple and undirected, an element of $E$, $(a,b)$,
 (or a 1-cell) can be
defined as the subset which contains two 0-cells $\{a,b\}$.  Intuitively,
a surface-cell (2-cell) represents a certain area, such as an unit square.
However, there is no ``area'' in a graph. So, a surface-cell can be defined
as a closed path of points, $p_0,...,p_n$, where $p_{i}$ and $p_{i+1}$ are adjacent
, or a path of edges, $e_0,...,e_n$, where $e_{i}$ and $e_{i+1}$ share a point.
A path forms a partial-graph of $G$. Thus, we have two choices:

First, if we define a surface-cell to be a closed path,
then the vertices of the path
may generate several closed paths. This approach can be used to
represent  and  describe indirect adjacency in $\Sigma_{m}$. For example
in $\Sigma_{2}$, assume $a,b, c,d $ are four points of an unit square,
then $\{a,b, c,d \}$ is a surface-cell. Moreover,  $\{a,b, c \}$, $\{a,b, d \}$, etc.
are (8-adjacency) surface-cells.  Such  cases may generate
ambiguities and local non-Jordan cases$^{6,7,8}$.

Second, a surface-cell must be not only a closed path (a partial graph),
but also a subgraph. That is
to say,  the simple closed path  has no any proper partial-graph in $G$ that is
a simple closed path. Because a subgraph is uniquely defined, we may
use the set of vertices of the subgraph to represent the surface-cell. This
matches our ideal, that is, any $i$-cells can be defined as a set of points.

[Note that  the original paper was confused with subgraph and partial graph. Shall be
swapped. ]

This paper mainly discusses the second case. i.e. a surface-cell is a subgraph
or a set of points in $G$. We will only consider this case
in the rest of the paper if we do not specify otherwise.
We sometimes call the first case the indirect-adjacency case as well as
the second case a direct-adjacency case $^{5,7,8}$.

According to the definition of subgraphs, if a subgraph $D$ of $G$ is
a simple cycle, then $D$ has no proper partial-graph
of itself that is a closed path. Such a $D$ is also
called a minimal-cycle. A surface-cell must be a minimal cycle,
but a minimal cycle may not be a surface-cell. Suppose that we have a $G$ as shown
in Figure 2.1.
In (a) and (b), every surface-cell is clear. However, $G$ in (c) and (d)
has the same $V$ and $E$, but each could have a different interpretation.
If the simple cycle
$(a,b,c,f,i,h,g,d)$ is not a surface-cell, then $G$ looks like a plane; if the
cycle is a surface-cell, then $G$ would have a 3-dimensional-cell (3-cell).
A geometrical interpretation of $G$ is needed to give a geometrical frame
(a topological structure) to $G$.

Let ${\cal C}$ be the set of all minimal cycles in $G$.
Defining a subset of ${\cal C}$ to be the set of surface-cells
is a way to generate a geometrical frame for $G$.

{\bf Definition 2.2}  Let ${\cal C}$ be the set of all minimal cycles in $G$.
A subset of ${\cal C}$, $U_{2}$ is called a surface-cell set if for any two different
minimal cycles in $U_{2}$, $u$ and $v$, $u \cap v$ is connected in $u \cap v$.

[The definition of the minimal cycle is still valid now. We want to add the following
definition to make the definition of surface-cell strong and more reasonable.
Each surface-cell or 2-cell is a minimal closed path.

Additional definition for discrete 1-cells and 2-cells: Each 1-cell contains a pseudo
vertex or point
in between two end vertices. This pseudo point is representative or reference point
of the 1-cell. Will not
allow any   edge to be added at the pseudo point.  A surface-cell also contains a pseudo
point. This pseudo point just like a central point for an object. It is used to refer that
this is a 2D or 1D entity or just a set of vertices.  It is reasonable since in Euclidean
plane, we always can make such a point. This idea was from our idea of refinement of
grid space. We can add finite number of points or vertices in these cells as necessary.
We can also define the those points in 3-cells and $m$-cells.  Only finite number of them
are allowed. They can be called the middle point, or inner points, or central point. They
are not real vertices. These points are called central pseudo points. These points are usually
used only indicating the status of cells and paths. Another type of  pseudo points
are relative pseudo points, for instance, a point only links two adjacent point, this point
can be deleted in the manifold. It is original point, not added later.
]

Each of the elements of $U_2$ is called surface-cell with respect to the pair
$<G,U_{2}>$. In the above definition, if $u \cap v$ is empty,
then $u$ and $v$ are not adjacent.
We want $u \cap v$ to be a connected path because we need the intersection of
two surface-cells to be  on the ``edges'' of these two surface-cells.
If $u \cap v$ is just a node (point), we say $u$ and $v$ are point-connected.
If $u \cap v$ has a line-cell, then   $u$ and $v$ are said to be
line-connected. See Fig. 2.1 (a) (b).

\begin{center}
\setlength{\unitlength}{0.01in}%
\begin{picture}(420,360)(40,460)
\thicklines
\put( 80,800){\line( 0,-1){ 40}}
\put(120,800){\line( 0, 1){  0}}
\put(120,800){\line( 0,-1){ 40}}
\put( 40,800){\line( 1, 0){160}}
\put(200,800){\line( 0,-1){ 40}}
\put(200,760){\line(-1, 0){160}}
\put( 40,760){\line( 0, 1){ 40}}
\put(160,800){\line( 0,-1){ 40}}
\put(280,780){\line( 0,-1){ 40}}
\put(280,740){\line( 1, 0){160}}
\put(440,740){\line( 0, 1){ 40}}
\put(440,780){\line(-1, 0){160}}
\put(280,780){\line( 1, 1){ 20}}
\put(300,800){\line( 1, 0){160}}
\put(460,800){\line( 0,-1){ 40}}
\put(460,760){\line(-1,-1){ 20}}
\put(440,780){\line( 1, 1){ 20}}
\put(320,780){\line( 1, 1){ 20}}
\put(360,780){\line( 1, 1){ 20}}
\put(400,780){\line( 1, 1){ 20}}
\put(320,740){\line( 0, 1){ 40}}
\put(360,740){\line( 0, 1){ 40}}
\put(400,740){\line( 0, 1){ 40}}
\put(320,740){\line( 1, 1){ 20}}
\multiput(340,760)(0.00000,7.27273){6}{\line( 0, 1){  3.636}}
\put(360,740){\line( 1, 1){ 20}}
\multiput(380,760)(0.00000,7.27273){6}{\line( 0, 1){  3.636}}
\put(400,740){\line( 1, 1){ 20}}
\multiput(420,760)(0.00000,7.27273){6}{\line( 0, 1){  3.636}}
\put(280,740){\line( 1, 1){ 20}}
\multiput(300,760)(0.00000,7.27273){6}{\line( 0, 1){  3.636}}
\multiput(300,760)(7.80488,0.00000){21}{\line( 1, 0){  3.902}}
\put(140,500){\framebox(40,0){}}
\put(120,500){\framebox(20,0){}}
\put(120,500){\framebox(60,0){}}
\put(140,620){\line( 0,-1){120}}
\put( 60,500){\framebox(120,120){}}
\put( 60,540){\line( 1, 0){120}}
\put( 60,580){\line( 1, 0){120}}
\put(100,620){\line( 0,-1){120}}
\put(298,500){\framebox(120,120){}}
\put(377,622){\line( 0,-1){120}}
\put(294,541){\line( 1, 0){120}}
\put(335,619){\line( 1,-3){ 12.300}}
\put(347,582){\line(-1,-5){  8.231}}
\put(300,576){\line( 6, 1){ 48.649}}
\put(349,582){\line( 5,-1){ 28.077}}
\put(378,578){\line( 1, 0){ 41}}
\put(338,540){\line( 0,-1){ 40}}
\put( 80,745){\makebox(0,0)[lb]{\raisebox{0pt}[0pt][0pt]{a}}}
\put( 80,805){\makebox(0,0)[lb]{\raisebox{0pt}[0pt][0pt]{b}}}
\put(120,805){\makebox(0,0)[lb]{\raisebox{0pt}[0pt][0pt]{c}}}
\put(120,745){\makebox(0,0)[lb]{\raisebox{0pt}[0pt][0pt]{d}}}
\put(315,785){\makebox(0,0)[lb]{\raisebox{0pt}[0pt][0pt]{b}}}
\put(315,725){\makebox(0,0)[lb]{\raisebox{0pt}[0pt][0pt]{a}}}
\put(340,805){\makebox(0,0)[lb]{\raisebox{0pt}[0pt][0pt]{c}}}
\put(345,765){\makebox(0,0)[lb]{\raisebox{0pt}[0pt][0pt]{d}}}
\put( 85,705){\makebox(0,0)[lb]{\raisebox{0pt}[0pt][0pt]{(a)}}}
\put(355,705){\makebox(0,0)[lb]{\raisebox{0pt}[0pt][0pt]{(b)}}}
\put(100,460){\makebox(0,0)[lb]{\raisebox{0pt}[0pt][0pt]{(c)}}}
\put(100,630){\makebox(0,0)[lb]{\raisebox{0pt}[0pt][0pt]{b}}}
\put(140,630){\makebox(0,0)[lb]{\raisebox{0pt}[0pt][0pt]{c}}}
\put( 60,630){\makebox(0,0)[lb]{\raisebox{0pt}[0pt][0pt]{a}}}
\put( 50,582){\makebox(0,0)[lb]{\raisebox{0pt}[0pt][0pt]{d}}}
\put( 89,581){\makebox(0,0)[lb]{\raisebox{0pt}[0pt][0pt]{e}}}
\put(131,581){\makebox(0,0)[lb]{\raisebox{0pt}[0pt][0pt]{f}}}
\put( 47,540){\makebox(0,0)[lb]{\raisebox{0pt}[0pt][0pt]{g}}}
\put(131,541){\makebox(0,0)[lb]{\raisebox{0pt}[0pt][0pt]{l}}}
\put( 91,542){\makebox(0,0)[lb]{\raisebox{0pt}[0pt][0pt]{h}}}
\put(292,629){\makebox(0,0)[lb]{\raisebox{0pt}[0pt][0pt]{a}}}
\put(336,626){\makebox(0,0)[lb]{\raisebox{0pt}[0pt][0pt]{b}}}
\put(374,630){\makebox(0,0)[lb]{\raisebox{0pt}[0pt][0pt]{c}}}
\put(280,581){\makebox(0,0)[lb]{\raisebox{0pt}[0pt][0pt]{d}}}
\put(331,583){\makebox(0,0)[lb]{\raisebox{0pt}[0pt][0pt]{e}}}
\put(368,581){\makebox(0,0)[lb]{\raisebox{0pt}[0pt][0pt]{f}}}
\put(281,542){\makebox(0,0)[lb]{\raisebox{0pt}[0pt][0pt]{g}}}
\put(330,546){\makebox(0,0)[lb]{\raisebox{0pt}[0pt][0pt]{h}}}
\put(370,543){\makebox(0,0)[lb]{\raisebox{0pt}[0pt][0pt]{l}}}
\put(334,461){\makebox(0,0)[lb]{\raisebox{0pt}[0pt][0pt]{(d)}}}
\end{picture}

\end{center}
\centerline{Figure 2.1 Examples of discrete spaces; (a) a minimal cycle in $2D$;}

\centerline{(b) a minimal cycle in $3D$; (c) minimal cycle $\{abcflhgda\}$
is not considered as a 2-cell;}

\centerline{(d) minimal cycle $\{abcflhgda\}$  is considered as a 2-cell.}

\vskip 0.2in

{\bf 2.2 Discrete Curves and Surfaces}

A subgraph is unique for a certain subset of vertices of $G$.
If $G'$ is a partial-graph, $G(G')$ is the subgraph with all
vertices in $G'$. For a subset of vertices of $G$, $V'$,  denote
by $G(V')$ the subgraph of $G$ with all vertices in $V'$.

{\bf Definition 2.2 } A semi-curve $D$ is a simple path $P=\{p_0,...,p_n\}$
such that $G(P)=P$ if $(p_0, p_n)$ is not an edge, or a semi-curve $D$ is
$G(P)=P\cup \{(p_0, p_n)\}$ if $(p_0, p_n)$ is an edge in
$G$.

We can view a semi-curve as a subset of vertices of $G$.
It is true that a semi-curve is a subgraph $D$ of $G$ where each vertex has one or
two adjacent vertices in $D$.  In other words,
a semi-curve $D$ is a simple path $p_0,...,p_n$ such that $p_i$ and $p_j$ are
not adjacent in $G$ if $i \ne j\pm 1 $ excepting $i=0$ and $j=n$.

{\bf Definition 2.3 }  For a graph $G=(V,E)$ and a $U_2$ of $G$,
 $D \subset V$ is called a curve if $D$ is a semi-curve and $D$ does not
 contain any surface-cell in $U_2$.

Note, we intentionally did not use $G(D)$ in Definition 2.3 because they are
equivalent.  We can summarize the above definitions to be one for the
discrete curves:

For a graph $G=(V,E)$ and a $U_2$ of $G$,
 $D \subset V$ is said to be a discrete curve with respect to $<G,U_{2}>$ if
(1) $D$ is connected (0-connected, or point-connected),
(2) each point (0-cell) is contained by one or two line-cells, and
(3) $D$ does not contain any surface-cell.

It follows easily that:

{\bf Lemma 2.1} If a semi-curve $C$ is not a curve, then $C$ is a surface-cell.

A common definition of curves is  that ``A simple discrete (or digital) curve is
just a simple path in $G$.''$^{9}$  It is a semi-curve in this paper. The reason is
that there is a simple path, such as $\{e,f,l,h\}$ in Fig. 2.1 (c), that cannot
separate the surface into two disconnected nonempty sets. This property is
called the Jordan theorem, which holds on continuous space. We will
discuss this property in section 3.

To determine the geometrical and topological structure of a graph $G$,
selection of $U_2$ is a critical issue.
In order to discuss the properties of all discrete spaces,
we give a default definition of $U_2$ here. Unless specified otherwise,
this default definition will apply.

{\bf Definition 2.4 } Let $G$ be a simple graph, the default definition
of $U_2$ of $G$ is :

(1) Assume $m$ is the minimum value of the lengths for all simple cycle.
	Then all simple cycle with length $m$ are surface-cells, i.e. in $U_2$.

(2) If a point $p$ is not a point of any surface-cell in $U_2$ by (1), a minimum
value of the simple cycle containing $p$ will be included in $U_{2}$.

After defining discrete curves, we now consider discrete surfaces.

{\bf Definition 2.5} Subgraph $D$ of $G$ is a semi-surface
if and only if each line-cell of $D$ is included in one or two
surface-cells in $D$.  $D$ is called a closed semi-surface
if and only if each line-cell is included in exact two surface-cells of $D$.

We also can view $D$ in Definition 2.5 as the vertex set of $D$ because
$D$ is uniquely defined by its vertex set. In other words, we can use
``vertex subset $D$'' substitutes ``subgraph $D$'' in the definition.
A 3-dimensional-cell (solid-cell or 3-cell) is a closed semi-surface,
but  a closed semi-surface
may not be a 3-cell. Similarly, we can define discrete surfaces
as follows. A closed semi-surface is called a minimal closed semi-surface
if it does not contain any proper subset that is a closed semi-surface.

{\bf Definition 2.6}  Let ${\cal S}$ be the set of all minimal closed semi-surfaces
in $G$. A subset of ${\cal S}$, $U_{3}$, is a 3-cell set if for
any two different elements in $U_{3}$, $u$ and $v$, $u \cap v$ is
empty, a vertex, point-connected if it contains line-cells, or
line-connected if it contains surface-cells (in  $u \cap v$).

In most cases, in  $u \cap v$ is a point (0-cell), a line-cell, or a surface-cell.
Each of the elements of $U_{3}$ is a 3-cell with respect to
$<G, U_{2}, U_{3}>$.


{\bf Definition 2.7}   Subset $S$ of $G$ is a discrete surface with  respect to
$<G, U_{2}, U_{3}>$ if and only if $S$ is a semi-surface and $S$ does not
contain any subset that is a 3-cell.

{\bf Definition 2.8} The boundary of surface $S$, denoted by $\partial S$, is
a subset of $S$ such that if  $b \in \partial S$ , then there is a line-cell
$B$ containing $b$ contained by  exactly
one surface-cell in $S$.

{\bf Corollary 2.1} A discrete surface is closed if and only if
$\partial S=\phi$. \\

{\bf 2.3 Regular Surface Points}

If $S$ is a subset of $G$, suppose that $S(p)$ contains all points in the
set of all surface-cells in $S$ containing $p$.

{\bf Definition 2.9}$^{10,11}$ A point $p$ in a discrete surface $S$ is regular
if the set of all surface-cells containing $p$ are line-connected among
these surface-cells.   If a point in $S$ is not regular, it is  called irregular.

We may generalize  the above definition. For any
$<G=(V,E),U_{2},U_{3}>$, a point $p\in V$ is said to be a regular surface
point if $S(p)$ (meaning the subgraph generated by $S(p)$ with all
line and surface-cells) is a surface and all surface-cells in $S(p)$ are
line-connected.

{\bf Lemma 2.5} For discrete surface $S$, let a point $p\in S$,
if $p$ has only two adjacent points $p',p''$ in $S$, then there are
two surface-cells $A,B$ such that $A\cap B$ contains  $p', p , p''$.
If $p',p''$ are adjacent in $S$, then  $p', p , p''$ form a surface-cell.

{\bf Lemma 2.6} Let $S$ be a  discrete surface. If $p$ is a inner and regular
point of $S$, then there exists a simple cycle containing
all points in $S(p)-\{p\}$ in $S$.

{\bf Definition 2.10} Let $S$ be a discrete surface. Then $p$ is said to be a
simple surface point if $S(p)-p$ is a closed curve.

A  regular inner point is a simple surface point in digital space $\Sigma_{m}$.
For a digital surface $S$ in $\Sigma_{3}$, a simple surface point under
Morgenthaler-Rosenfeld's definition  is a regular inner surface point in the case of
direct adjacency$^{8}$. \\\\

\section{ Simply Connected Discrete Surfaces and the Jordan Curve Theorem}

In topology, there is a fundamental theorem called the Jordan curve theorem:
A simply closed curve $J$ in a plane $\Pi$ decomposes $\Pi - J$ into two
components$^{12}$.  In fact, this theorem holds for a simply connected surface.
A plane is a simply connected surface in continuous space. This theorem is not true
for a general continuous surface. For example, a meridian circle of a ring-surface
can not separate the ring-surface into two components.

In continuous space, a connected topological space $T$ is simply connected
if for any point $p$
in $T$, any simply closed curve containing $p$ can be contracted to $p$.
The contraction is a continuous mapping among a series of closed curves.
To simulate the proof of the Jordan theorem in discrete space is not easy. First,
we must define a ``discretely'' continuous mapping, then we need the concept
of ``discrete contraction.''\\

{\bf 3.1 Gradual Variation and Simply Connected Surfaces}

Rosenfeld proposed ``continuous'' functions on digital spaces.$^{13}$
Pawlak used the concept of ``roughly continuous.''$^{14}$
Chen proposed $\lambda$-connectedness to represent the
relationship between a pixel and its adjacent pixel in an image. $^{15}$
In order to simplify the proof of an important interpolation and extension
theorem, Chen translated $\lambda$-connectedness to so called
``gradual variation.''$^{16}$ In this paper, we still use gradual variation
to represent the ``continuous'' contraction. On the otherhand, a simply
connected continuous
surface is orientable. In this section, we assume the discrete surface is
both regular and orientable. $^{5}$

{\bf Definition 3.1}$^{17}$ Let $G$ and $G'$ be two connected graphs. A mapping
$f: G\rightarrow G'$ is gradually varied if for two vertices
$a, b \in G$ that are adjacent in $G$,
then $f(a)$ and $f(b)$ are adjacent in $G'$ or $f(a) = f(a')$.

In this paper, we do not directly use the gradually varied function, but
we need to use the same idea to build a concept for ``continuously'' changing
from a curve to another curve.  ``Gradual
variation'' is still a good term to describe the ``continuous''
change between two simple paths. Herman defined ``elementarily N-equivalent''
for defining simply connected space.$^{1}$

We will call a simple path a pseudo-curve in the rest of the paper in
order to extend the concept to 2D and higher dimensional objects.
Intuitively, ``continuous'' change from a simple path $C$ to another $C'$ is that
there is no ``jump'' between these two paths.  If $x, y\in S$, $d(x,y)$ denotes
the distance between  $x$ and $y$. $d(x,y)=1$ means that $x$ and $y$ are
adjacent in $S$.

{\bf Definition 3.2} Two simple paths $C=p_0,...,p_n$ and
$C'=q_0,...,q_m$ are gradually varied
in (a regular surface) $S$ if  $d(p_0,q_0)\le 1 $ and  $d(p_n,q_m)\le 1 $
and  for any non-end point $p$ in $C$, then

(1) $p$ is in $C'$, or $p$ is contained by a surface-cell $A$ (in $G(C\cup C')$)
	  such that $A$ has a point in $C'$.

(2) Each non-end-edge in $C$ is contained by a surface-cell $A$ (in $G(C\cup C')$) which
	 has an edge contained by $C'$ but not $C$ if $C'$ is not a single point.

\noindent and vise versa for $C'$.

For example, $C$ and $C'$ in Fig. 3.1 (a) are gradually varied,
but $C$ and $C'$ in Fig. 3.1 (b) are not gradually varied.
We can see that a surface-cell, which is a simple path, and any point in the
surface-cell are gradually varied.  Assume $E(C)$ denotes all edges in
path $C$.  Let $XorSum(C,C') = (E(C)-E(C'))\cup (E(C')-E(C))$.
$XorSum$  is called $sum ( modulo 2)$ in  Newman's book $^{18}$.

\begin{figure}[h]
	\begin{center}

   \epsfxsize=3.5in
   \epsfbox{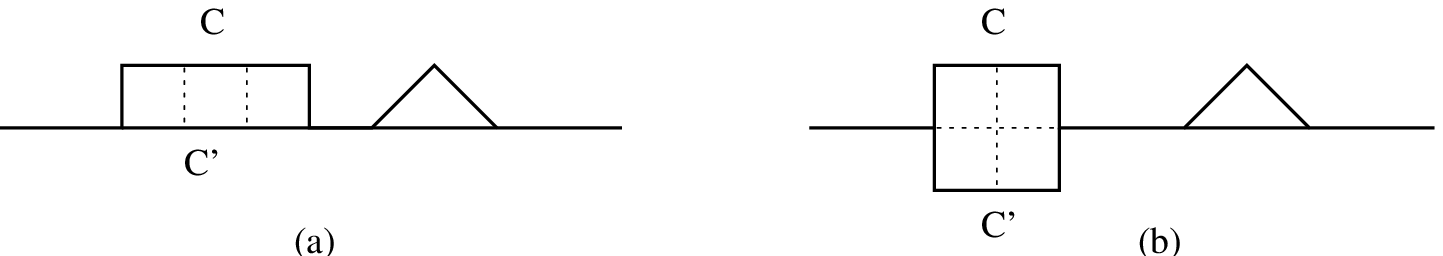 }

\caption*{Figure 3.1 (a) $C$ and $C'$ are gradually varied; (b) $C$ and $C'$ are not gradually varied.}

	\end{center}
\end{figure}



{\bf Lemma 3.1} Let $C$ be a  pseudo-curve and $A$ be a surface-cell.
If $A \cap C$ is an arc containing at least an edge,
then $XorSum(C,A)$ is a gradual variation of $C$.

It is not difficult to see that $XorSum(XorSum(C,A),A)= C$ and
$XorSum(XorSum(C,A),C)= A$ under the condition of Lemma 3.1.

{\bf Definition 3.3} Two simple paths (or pseudo-curves) $C,C'$ are said to
be homotopic if
there is a series of simple paths $C_0,...,C_n$ such that
$C=C_0$, $C'=C_n$, and  $C_{i}, C_{i+1}$ are gradually varied.

{\bf Lemma 3.2} If two simple paths $C,C'$ are homotopic then
there is a series of simple paths $C_0,...,C_m$ such that
$C=C_0$, $C'=C_n$, and  $XorSum(C_{i}, C_{i+1})$ is a surface-cell excepting
end-edges of $C,C'$.

Because a surface-cell $A$ is a closed path, we can define two
orientations (normals ) to $A$: clockwise and counter-clockwise.
Usually, the orientation of a surface-cell is not a critical issue.
However, for the proof of the Jordan theorem it seems necessary. In fact, a curve
which is a set of points has no ``direction,'' but a pseudo-curve, a path,
has its own ``travel direction'' from $p_0$ to $p_n$.  For two paths $C$ and $C'$,
which are gradually varied, if a surface-cell $A$
is in $G(C\cup C')$, the orientation of $A$ with respect to $C$ is determined by
the first pair of points $(p,q) \in C \cap A$ and $C= ... p q ... $ .
Moreover, if a line-cell of $A$ is in $C$ (meaning all 1-cells of $A$ are in $C$ ),
then the orientation of $A$ is fixed with respect to $C$.

According to Lemma 2.6, $S(p)$ contains all adjacent points of $p$ and
$S(p)-\{p\}$ is a simple cycle
(there is a cycle containing all points in $S(p)-\{p\}$).
We assume that cycle $S(p)-\{p\}$ is always oriented
clockwise. For two points $a, b \in S(p)-\{p\}$,
there are two simple cycles containing the path $a\rightarrow p \rightarrow b$ :
(1) a cycle from $a$ to $p$ to $b$ then moving clockwise to $a$, and
(2) a cycle from $a$ to $p$ to $b$ then moving counter-clockwise to $a$.
See Fig. 3.2.  It is easy to see that the simple cycle $S(p)-\{p\}$ separates
$S-\{S(p)-\{p\}\}$ into at least two connected components because from
$p$ to any other points in $S$ the path must contain a point in $S(p)-\{p\}$.
$S(p)-\{p\}$ is called a Jordan curve.


\begin{figure}[h]
	\begin{center}

   \epsfbox{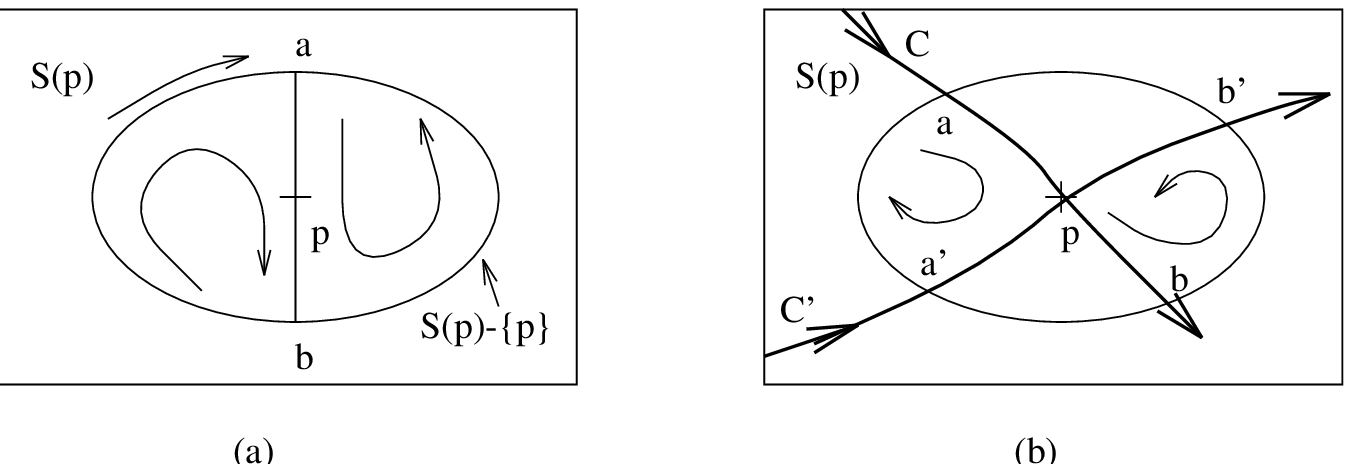 }

\caption*{Figure 3.2 (a) Two adjacent points $a$ and $b$ of $p$ in $S(p)$; (b) an example for two cross-over paths.}

\end{center}
\end{figure}

{\bf Definition 3.4} $C$ and $C'$ are said to  ``cross-over'' each other if
there are  points $p$ and $q$ ($p$ may be the same as $q$) such that
$C = ... a p b...s q t...$ and $C'= ...a' p b...s q t'...$ where
$a\ne a'$ and $t\ne t'$. The cycle
$a p a' ... a$ without $b$ in $S(p)$ and the cycle $q t ...t' q$ without
$s$ in $S(q)$ have different orientations with respect to $C$.

For example, in Fig. 3.2, $C$ and $C'$ are ``cross-over''
If  $C$ and $C'$ are not ``cross-over'' , then we say that $C$ is at a side of
$C'$.

{\bf Lemma 3.3}  If $C$ and $C'$ do not cross-over each other, and they
are gradually varied, then every surface-cell in $G(C\cup C')$ has the
same orientation with respect to the ``travel direction'' of $C$ and
opposite to the ``travel direction'' of $C'$.

We also say that $C,C'$ in Lemma 3.3 are sided-gradually varied.
Intuitively, a simply connected set is such that for any point, every cycle containing
the point can contract to the point.
According to the nature of the word ``contraction,'' we can give the mathematical
definition of "contraction" for discrete spaces. In fact, the contraction
procedure relates to some substance. This substance gradually loses the
size (which can be space occupied), and when a part was lost, it will
never come back again.

{\bf Definition 3.5}  A simple cycle $C$ can contract to a point $p\in C$
if there exist a series of simple cycle, $C=C_0, ..., p=C_{n}$:
(1) $C_i$ contains $p$ for all $i$;
(2) If $q$ is not in $C_i$ then $q$ is not in all $C_{j}$, $j>i$;
(3) $C_{i}$ and $C_{i+1}$ are side-gradually varied.

We now show three reasonable definitions of simply connected spaces below.
A general definition of a simply connected space should be :

{\bf Definition 3.6(a)} $<G,U_2>$ is
 simply connected if any two closed simple paths are homotopic.

However, if we use this definition, then we may need an extremely long proof for
the Jordan theorem. The next one is the special case of the Definition 3.6(a),

{\bf Definition 3.6(b)}  A connected discrete space $<G,U_2>$ is  simply connected
if for any point $p\in S$, every cycle containing $p$ can contract to $p$.

This definition of the simply connected set is based on the original meaning of
simple contraction. In order to
make the task of proving the Jordan theorem
simpler, we give the third strict definition of simply connected surfaces as
follows.

We know that a simple closed path (simple cycle) has at least
three vertices in a simple graph.
This is true for a discrete curve in a simply connected surface $S$.
For simplicity, we call an unclosed path an arc.
Assume $C$ is a simple cycle with clockwise orientation.
Let two distinct points $p, q\in C$. Let $C(p,q)$ be an arc of $C$ from $p$ to $q$ in
a clockwise direction, and $C(q,p)$ be the arc from $q$ to $p$ also in a clockwise
direction,
then we know $C=  C(p,q)\cup C(q,p)$. We use $C^{a}(p,q)$  to represent
the counter-clockwise arc from $p$ to $q$. Indeed, $C(p,q) = C^{a}(q,p)$.
We always assume that $C$ is in clockwise orientation.

{\bf Definition 3.6(c)} A connected discrete space $<G,U_2>$ is simply connected
if for any simple cycle $C$
and  two points $p, q\in C$, there exists a side-gradually varied
simple cycle path
$Q_{0},...,Q_{n}$ such that $C(p,q)=Q_{0}$ and $C^{a}(p,q)=Q_{n}$.

[In fact, Definition 3.6(b),(c) are special cases of  Definition 3.6(a). Definition 3.6(b) and
Definition 3.6 (c) are equivalent. In Chen's book 2004,  we think the proof of the equivalence
is hard. It is not. we can now extend a point as two points, so  Definition 3.6(b) is a special
case of Definition 3.6(c). ]

{\bf Proposition 3.1 (new)}  Definition 3.6(b) and Definition 3.6 (c) are equivalent.

{\bf Proof} Definition 3.6(b) is a special case of Definition 3.6(c).
Now we prove that Definition 3.6(c) can be induced from Definition 3.6(b).
When we select contracting point $x$ is  $p$  , then let the contracting
sequence $C_0$, $C_1$,...,$C_i$ contain $q$, (so path from $p$ to $q $, $C_{k}(p,q)$, and $q$ to $p$,  $C_{k}(q,p)$
have their own corresponding gradually varied paths, respectively in $C_k$, $k=0$,...,$i$.)
but $C_{i+1}$, ..., $C_{n}=p$ do not contain $q$. ($C_{i+1}$ does not contain $q$ is the key,
the other are not really matter.) Note that, all $C_{t}$, $t=i+1$,...,$n-1$ are closed path.
In $C_{i+1}$, $q$ has a corresponding point  in $C_{i+1}$, say $q^{i+1}$, ($q$ changed to  $q^{i+1}$, in the process.)
from $p$ to $q^{i+1}$ . There are two paths $p$ to $q^{i+1}$,  $C_{i+1}(p,q^{i+1})$, and $q^{i+1}$ to $p$,
$C_{i+1}(q^{i+1},p)$ . Therefore, $C_{i}(p,q)$ and $C_{i+1}(p,q^{i+1})$ are gradually varied, so are
$C_{i}(q,p)$ and $C_{i+1}(q^{i+1},p)$ .  In the same way, we can find $q^{i+2}$,...,$q^{n-1}$ . Thus,
$C_{0}(p,q)$,...,$C_{i+1}(p,q^{i+1})$,...,$C_{n-1}(p,q^{n-1})$, $C_{n-1}(q^{n-1},p)$,... $C_{i+1}(q^{i+1},p)$ ,$C_{0}(q,p)$
are such a gradually varied sequence.

Basically, the deformation does not really care about cross-over points. Does not allow cross-over points will
make the proof easier.

[More discussion : If we can prove the follow statement then Definition 3.6(b) and
Definition 3.6(c) are equivalent.

`` For a simply connected surface by Definition 3.6(b), let $C$ be a simple cycle
and let two points $p, q\in C$.
Then there is a simple cycle path,
$C(p,q)=Q_{0},...,C^{a}(p,q)=Q_{n}$, such that $Q_{i}$ and $Q_{i+1}$
are side-gradually varied for all $i$. '']\\

{\bf  3.2 The Jordan Theorem }

Since a simple cycle could be a surface-cell, it can not separate
$S$ into two disconnected components. [In the case of allowing the central pseudo points,
we will have the general Jordan Curve Theorem. We will prove that in the last of this
section]

However, for a closed discrete curve, we have

{\bf Theorem 3.1} (The Jordan Theorem) Discrete simply connected surfaces $S$,
defined by Definition 3.6(c),
have Jordan's properties:
A closed discrete curve $C$ which does not contain any point of $\partial S$
divides $S$ into at least two disconnected components.
In other words, $S-C$ consists of at least two disconnected components.
(These components are disconnected.)

{\bf Proof}  Suppose that $C$ is a closed curve in a simply connected
surface $S$. $C$ does not reach the border of $S$, i.e.
$C\cap {\partial S} = \emptyset$.
 Assume point $p\in C$, then suppose that $q$ and $r$
are two adjacent points of $p$ in $C$ with form of $...q p r,...$, where
the direction of  ... $q$ to $p$ to $r$ ...to $p$ is clockwise.
See Fig. 3.3.
$\{p,r\}$ is a line-cell, then there are two surface-cells containing
$\{p,r\}$. Denote these by $A$ and $B$ with clockwise orientation.
Our strategy is to prove that if
there is a point $a$ in $A$ which is not in $C$, and a point $b\in B$
and $b\notin C$, then any path from $a$ to $b$ must contain a point
in $C$. Then we can see that $S-C$ are not (point-) connected and we have the
Jordan theorem.

First, we want to prove that there must exist a point in $A-C$.
 If each point in $A$ is in $C$, since $A$ is a simple cycle,
then $C=A$. However, $C$ is not a surface-cell, so the statement can
not be true. Thus, there is a point $a\in A-C$. For the same reason
there is a point $b\in B-C$.   We assume
that $a$ is the last such point in $A$ starting with $p$, and $b$ is
the first such  point in $B$ starting with $p$. (see Fig 3.3) We always
assume  clockwise direction here unless we indicate otherwise.

On the other hand, based on the discussion of Section 3,
if a surface $S$ does not contain any abundant point then the intersection of
any two surface-cells in $S$ at most contains a line-cell. So, $A\cap B = \{q,p\}$.
There is a path from $a$ to $b$, $Q(a,b)=a p q...q' b$, where
$p, q,...,q'$ are in $C$.

Suppose we make the counter statement: there is a path from $a$ to $b$, $P(a,b)$,
such that there is no point of $P(a,b)$ in $C$. Because  $Q(a,b)$ has
only two points $a,b$ in $S-C$,  $P(a,b)\cap Q(a,b)$
just contains two points $\{a,b\}$. Thus,  $D=P(a,b)\cup Q^{(a)}(b,a)$ is a simple
cycle. We may assume $P(a,b)\cup Q^{(a)}(b,a)$ is oriented clockwise.

According to Definition 3.6(c), For simple cycle $D$,  there are
finite simple cycle $D=D_0,...,D_{m}=\{p\}$ where $D_{i}$ and
$D_{i+1}$ are side-gradually varied. (We suppose $D_{i}$ and $D_{j}$ are
different). Meanwhile, there are
finite paths $P(a,b)=P_{0}(a,b),..., P_{n}(a,b)=Q(a,b)$ so that
$P_{i}(a,b)$ and $P_{i+1}(a,b)$ are side-gradually varied.


\begin{figure}[h]
	\begin{center}

   \epsfbox{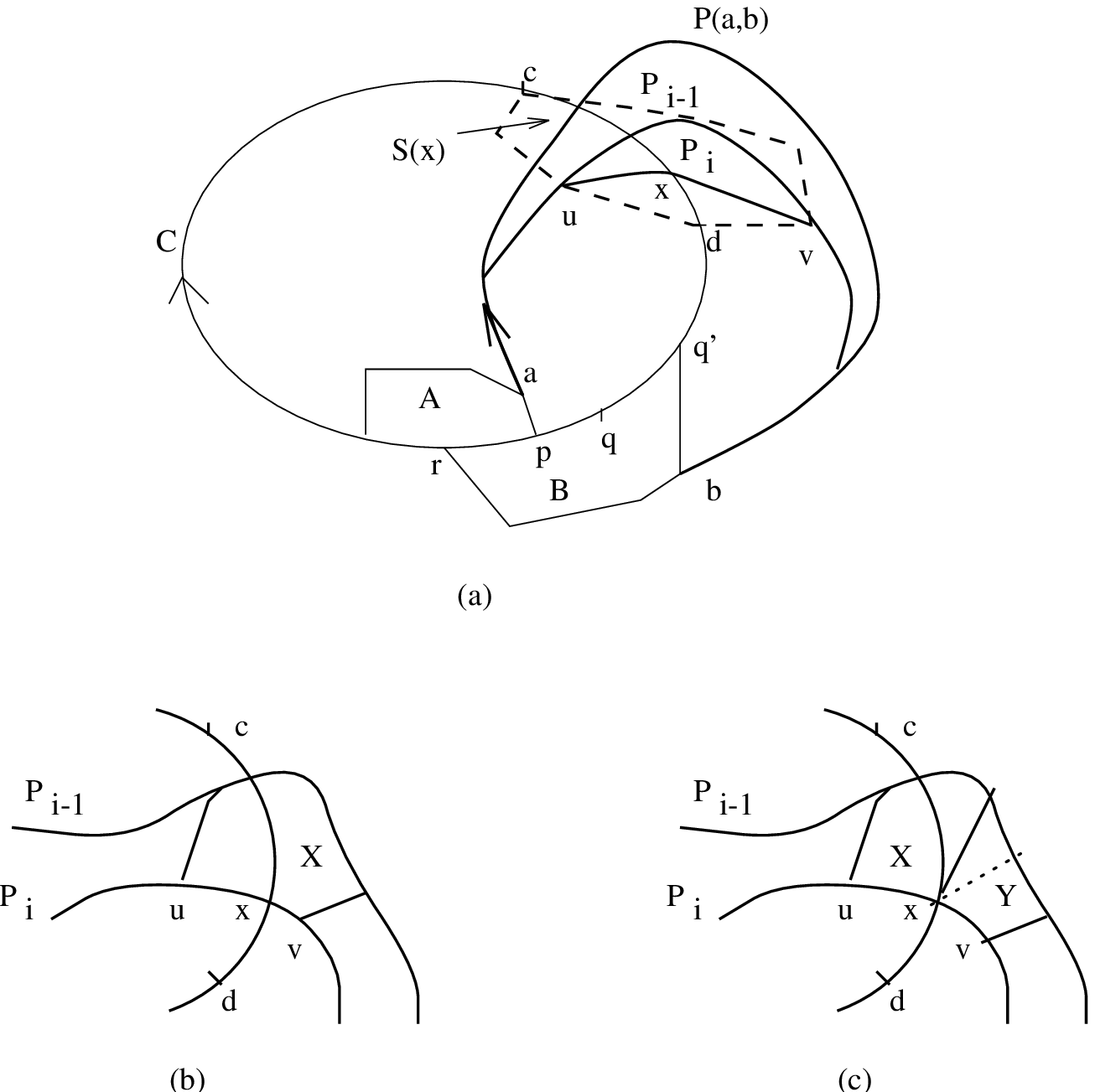 }

\caption*{Figure 3.3 A close curve $C$ and the paths from $a$ to $b$}
\end{center}
\end{figure}

Because  $P_{n}(a,b)=Q(a,b)$ is cross-over to $C$,
we only need to prove the following statement to reach a contradiction:
each $P_{j}(a,b)$ does not cross-over $C$.

If there is a $P_{i}(a,b)$ cross-over to $C$, let $P_{i}(a,b)$ be the first
one, i.e.,  $P_{i-1}(a,b)$ do not cross-over $C$.
Let point $x$ in $P_{i}(a,b) \cap C$ and $x\notin P_{i-1}(a,b)$.

Suppose that $x=``p''$ and $x=``q''$ in Definition 3.4 and assume
$P_{i}(a,b)=...u x v...$ and $C= ... c x d ...$, where $v\ne d$. We know
that $u,v,c,d $ are in simple cycle $S(x)-\{x\}$ (Lemma 2.6).

According to Definition 3.2, there is a surface-cell $X$ contains $(u,x)$ and
an edge $e\in P_{i-1}(a,b)-P_{i}(a,b)$. We might as well assume that $x$ is the first
point on $P_{i}(a,b)$ that is in $C$. Thus, $c,d\notin X$. If $X$ also contains
$v$, then $c,d$ in $S(x)$ must in the same side of $u x v$. Therefore,
$C$ and $P_{i}(a,b)$
do not cross-over each other at $x$. (See Fig. 3.3 (b).) If $X$ does not contain
$v$, then there is a 2-cell $Y$
contains $(x,v)$. We can see  that $X$ and $Y$ are line-connected in $S(x)$;
otherwise, there is an edge in $P_{i-1}(a,b)$  contained by a 2-cell in
$G(P_{i-1}(a,b) \cup  P_{i}(a,b))$ which  has no edge in $P_{i}(a,b)$.
(See Fig. 3.3 (c).) Suppose the $X\cap Y= (x,y)$, then $y$ is on $P_{i-1}(a,b)$
and $u...y...v$ is on the simple cycle $S(x)-\{x\}$.
$c,d$ in $S(x)$ must in the same side of $u x v$. Therefore, $C$ and $P_{i}(a,b)$
do not cross-over each other at $x$.

If $C=...c x s...t y d... $ and $P_{i}(a,b)= ...u x s...t q v$
where $x\ne y$, $c\ne u$, and $d\ne v$. We know that  none of $x, s,...,t, y$
is in $P_{i-1}(a,b)$. Each surface-cell in $G(P_{i-1}(a,b) \cup  P_{i}(a,b))$
containing $(u,x)$, $(x,s)$,..., or
$(y,t)$ does not contain $c$ or $d$. $C$ is on the side of
$G(P_{i-1}(a,b) \cup  P_{i}(a,b))$ at the arc of $c x s...t y d$.
Consider $S(x)$ and $S(y)$, $C$ does not cross-over $P_{i}(a,b)$.

However, $P(n)$ do cross-over $C$. Thus, there is no path
from $a$ to $b$ that does not contain a point of $C$.
We have proved of Jordan theorem for simply connected discrete
surfaces. $\diamondsuit$\\

{\bf The New and More Detailed Proof:}
Here is the new modified proof.   Part 1 is just clarify some statements. Part 2 has
more changes.

Suppose that $C$ is a closed curve (or simple path, at least
has a pseudo point in the center of a 2-cell)

The idea of the proof is to find two points in each sides of curve $C$.
This is because that for any 1-cell $(r,p)$ in $C$, there are two 2-cells $A$,$B$
sharing $(r,p)$ by 2D discrete surface definition.  $A$ must contain a vertex $a$ and $B$ must contain $b$,
and they are not in $C$. ($a$, $b$ are adjacent to some point in $C$).
We are going to prove that from $a$ to $b$, a path must cross-over $C$.
That is the most important part of the Jordan curve theorem.

We assume, on the contrary, there is a simple path from $a$ to $b$ does not cross-over $C$, called
$P_(a,b)$. But we know there is $P(b,a)$ in $A\cup B$ does cross-over $C$. (Fig. 3.3(a))

Because $S(r)$ containing all 2-cells that contains
$r$, the boundary of $S(r)$ is a simple curve. This is because
 we always assume that $r$ is a regular point. (When $a$ is not a pseudo-point,) $a$ is on the boundary of $S(r)$.
(The boundary of $S(r)$ is denoted as $S(r)-{r}$). $A\cup B$ is a subset of $S(r)$.

We can assume that $P_(a,b)\notin S(r)$; otherwise, it must cross-over $C$. (a 2-cell containing $r$ must
have an edge on $C$, or all points of the 2-cell are on the boundary of $S(r)$ except $r$). If $P(a,b)$ does
not contain $r$, must be a part of boundary of $S(r)$ which is a cycle.  $r$ has two adjacent points on $C$, (if
they are not pseudo, meaning here it can be eliminated or added on an edge
 that does not affect to the 2-cell) so these two points are also in the boundary of  $S(r)$.
So there are only two ways from a to be on the boundary of $S(r)$. These two points are not on the same side of the cross-over path containing $r$. (The boundary of $S(r)$ was separated by
  the cross-over path containing $r$.) $P_(a,b)$ must contain a such point that is on $C$.

Therefore, $P_(a,b)\notin S(r)$. Then $P_(a,b)\cup P(b,a)$ is a simple closed curve. ($P(b,a)$ passes $r$).
By the definition of the simply-connected surface,
there are
finite numbers of paths $P(a,b)=P_{0}(a,b)$,..., $P_{n-1}(a,b)$, such that
so that $P_{i}(a,b)$ and $P_{i+1}(a,b)$ are (side-)gradually varied.
In addition $P_{n-1}(a,b)$ is gradually varied
to $P_{n}(a,b)=P^{a}(b,a)$ (reversed $P(b,a)$ that passes $r$).

We now can assume that there is a smallest $i$ such that $P_{i}(a,b)$ cross over $C$, but $P_{i-1}(a,b)$ does not.
(Fig. 3.3 (a)). The idea is we will prove that it is impossible if $P_{i-1}(a,b)$ does not cross over $C$.

Let point $x$ in $P_{i}(a,b) \cap C$ and $x\notin P_{i-1}(a,b)$. There are two cases:
(1) cross over single point, or (2) cross over a sequence of points on $C$.

{\bf Case 1}: Suppose that $x=``p''$ and $x=``q''$ in Definition 3.4 and assume
$P_{i}(a,b)=...u x v...$ and $C= ... c x d ...$, where $v\ne d$.

We know that $u,v,c,d $ are in the boundary of $S(x)$, a simple cycle $S(x)-\{x\}$ (Lemma 2.6).
There is a surface-cell $X$ (in between $P_{i-1}$ and $P_{i}$) contains $(u,x)$.
$X$ has a sequence of points $S1$ in $P_{i-1}$ and  a sequence of points $S2$ in $P_{i}$.
$X$ has at most two edges $e1$, $e2$ not in  $P_{i-1}\cup P_{i}$ ; $S1$, $e1$, $S2$, $e2$,
are the boundary of $X$. $e1$ is the edge linking $S1$ to $S2$,  and $e2$ is the edge linking $S2$ to $S1$
counterclockwise.(Again, $e1$ may or may not be directly incident to $u$, and $e1$ may be an empty edge if
$P_{i-1}$ intersects $P_{i}$ at point $u$. $e2$ may also in the same situation.)
We might as well assume that $x$ is the first
point on $P_{i}(a,b)$ (from $a$ to $b$ in path $P_{i}$ )that is in $C$. Thus, $c,d\notin X$.
(If $c$ is in $X$ $c$ must be in $P_{i-1}$. if $d$ is in $X$, $x$ is not only cross over point. )

If $X$ contains $v$,  we will have a cycle $u\cdot d \cdot v (e2) (S1) (e1) $ in the boundary of $S(x)$
$(e2)(S1) (e1)$ contains only points in $P_{i-1}$ and $u$,$v$ (that are possible end points of $e1$, $e2$).
$c$ is on the boundary of $S(x)$ too. Where is $c$? It must be in the boundary curves (of $S(x)$)
from $u$ to $d$ or the curve from $d$ to $v$.
Then $c,d$ in $S(x)$ must in the same side of $u x v$ which is part of $P_{i}$.
Therefore, $C$ and $P_{i}(a,b)$ do not cross-over each other at $x$. (See Fig. 3.3 (b).)

If $X$ does not contain $v$, then there must be a 2-cell $Y$ (in between $P_{i-1}$ and $P_{i}$)
containing $(x,v)$. We can see that $X$ and $Y$ are line-connected in $S(x)$.
(See Fig. 3.3 (c).)
This is due to the definition of regular point of $x$, all surface-cells containing
$x$ are line-connected. Meaning there is a 2-cell paths they share a 1-cell
in adjacent pairs.

Since $X$ and $Y$ are line-connected, we can assume: a) $X\cap Y= (x,y)$, then $y$ is on $P_{i-1}(a,b)$.
Let $e3$ be the possible edge from $v$ to $P_{i-1}(a,b)$. ($e3$ could be empty as $e1$)
and $u(e1)..y...(e2)v$ is on the boundary cycle of $S(x)$.  Except $u$ and $v$, $u(e1)..y...(e2)v$ is
on $P_{i-1}(a,b)$. $u...d...v$ is part of the boundary cycle of $S(x)$.
In addition, $c$ (that is not in $P_{i-1}(a,b)$) must be in the boundary curves (of $S(x)$)
from $u$ to $d$ or the curve from $d$ to $v$. Again,
$c,d$ in $S(x)$ must in the same side of $u x v$ which is part of $P_{i}$.
Therefore, $C$ and $P_{i}(a,b)$ do not cross-over each other at $x$. (See Fig. 3.3 (c).)
b) $X\cap Y= x$ , let us assume that $e1$ incident to $P_{i-1}(a,b)$ at $y'$ ($y'$ is $u$ if
 $e1$ is empty. ) and $e3$ incident to $P_{i-1}(a,b)$ at $y''$.
 We will have a set of points $y'=y_0, y_1,...,y_k=y''$ in $P_{i-1}(a,b)$. All $y_i$'s are contained
 in a 2-cell containing $x$. All $y_0, y_1,...,y_k$ are in the boundary cycle of $S(x)$.
 $c$ that is not in $P_{i-1}(a,b)$. $c$ must be in the boundary curves (of $S(x)$)
from $u$ to $d$ or the curve from $d$ to $v$. Thus,
$c,d$ in $S(x)$ must in the same side of $u x v$ which is part of $P_{i}$.
$C$ and $P_{i}(a,b)$ do not cross-over each other at $x$. (See Fig. 3.3 (c).)

{\bf Case 2:} Suppose $P_{i}(a,b)$ and $C$ cross over a sequence of points on $C$:
$P_{i}(a,b)=...u x_0 x_1...x_m v...$ and $C= ... c x_0 x_1...x_m d ...$, where $v\ne d$.

We still have $e1 = (u,y_0)$ and $e3=(v,y_k)$ where $y_0$ and $y_k$ are on $P_{n-1}$ for some $k$
all $y_t$ , $t=0,1,...,k $, are in a 2-cell that containing $x_j$, $j=0,1,...,m$.
(If $u$ does not have a direct edge linking to $P_{n-1}$, $u$ will be in a 2-cell between $P_n$ and
$P_{n-1}$,  either $u$ is a pseudo point on $P_{n}$ for the deformation from $P_{n-1}$ to $P_n$, or
 $P_{n-1}$ and $P_n$ intersects at $u$. That $u$ is a pseudo point means here it has a neighbor that
 has an edge link to $P_{n-1}$, or the neighbor's neighbor, and so on. We can just assume here
 $u$ is the point that is adjacent to a point in $P_{n-1}$. In the theory, as long as
 $u$ is contained by a 2-cell such that all the points in the 2-cell are in $P_{n-1}$ or $P_n$.)

The same way will apply to this case just treat $x_0$,...,$x_m$ to $x$ in Case 1.  We first get the union
of $S(x_0)$,...,$S(x_m)$.

The boundary of this union will be simple cycle too; using mathematical induction
we can prove it. Then, we can prove the rest of theorem using the same method in Case 1. See Fig. 3.4.

\begin{figure}[h]
	\begin{center}

   \epsfbox{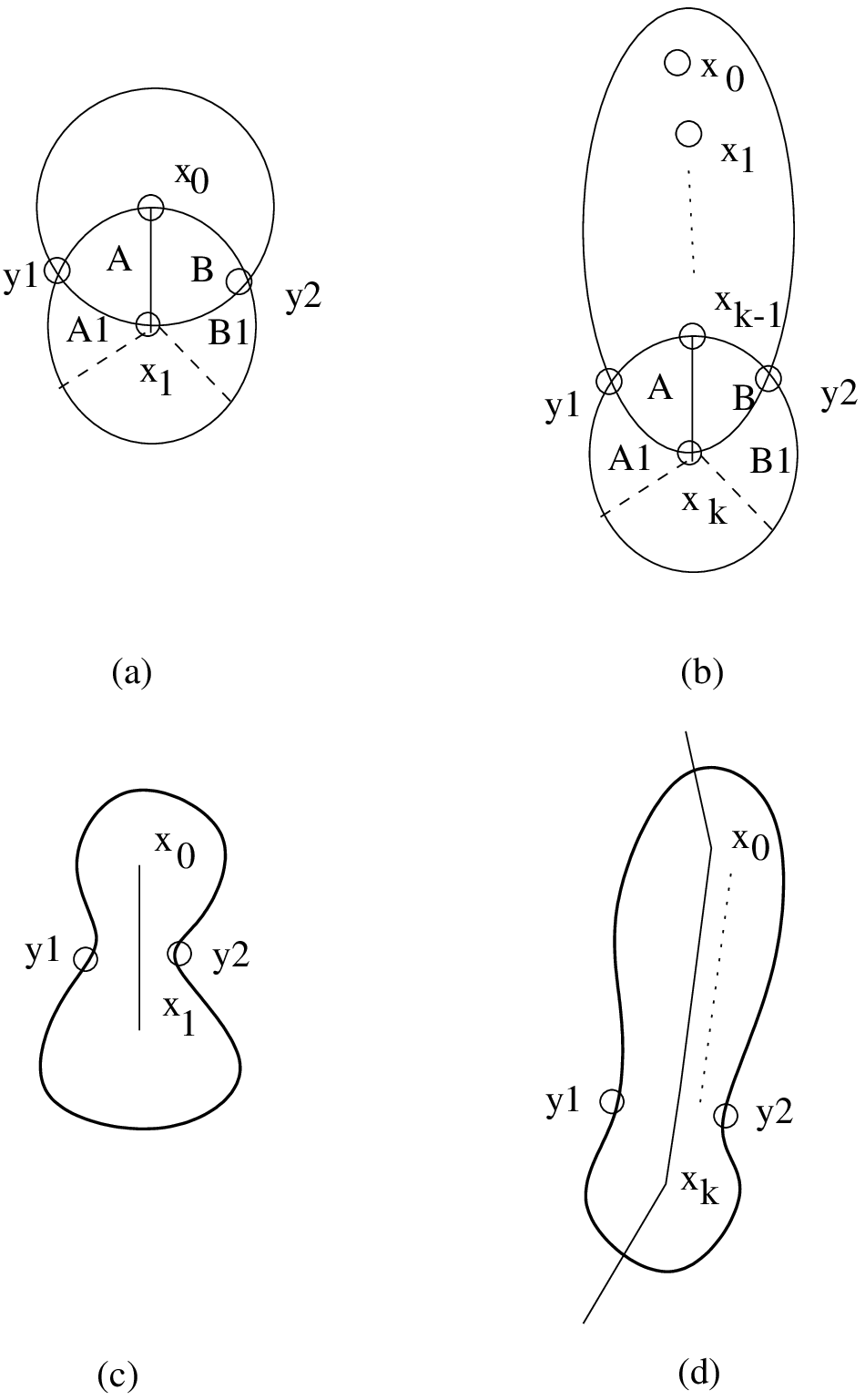 }

\caption*{Figure 3.4  The union of neighborhoods of a sequence of adjacent points, $S(x_0,...,x_k)$ and its boundary}
\end{center}
\end{figure}

The following is the detailed proof: (Some idea was represented in the original proof.)

Let $S(x_0,...,x_k)=S(x_0)\cup...\cup S(x_k)$. We will first prove that the boundary of $S(x_0)\cup S(x_1)$ is a simple
cycle (it is a simple closed curve too).

We know that $(x_0, x_1)$ is an edge in $C\cap P_{i}(a,b)$. Also, there are two 2-cells $A,B$ in $S(x_0)$ containing  $(x_0, x_1)$ .

$x_1$ is a boundary point in $S(x_0)$ , so no other 2-cell will contain $x_1$. In the same way, $S(x_1)$ also contains  $A,B$, and
$x_1$ is only contained in two 2-cells in $S(x_1)$. Therefore,  $S(x_0)\cap S(x_1)= A\cup B$ and  $A\cap B=(x_0,x_1)$.

$A$ and $B$ are adjacent 2-cells. On the other hand, $x_1$ is on the boundary curve (that is closed) of $S(x_0)$, so $x_1$ has two adjacent points
on this cycle, $y_1$ and $y_2$. (We assume that $y_1$ and $y_2$ are not pseudo points, so)  $y_1$ and $y_2$ are both on the boundary
of $S(x_0)\cup S(x_1)$. (If $y_1$ or $y_1$ is pseudo points, we can ignore $y_1$ or $y_2$ to find the a actual point that adjacent to $x_1$.)
$(x_1,y_1)$ has two 2-cells containing $(x_1,y_1)$ in $S(x_0)\cup S(x_1)$. For instance, in Fig. 3.4 (a) , $A$ and $A_1$ contain $(x_1,y_1)$ and
$ B$ and $B_1$ contain $(x_1,y_2)$. Thus, the boundary of $S(x_0)\cup S(x_1)$ is a closed curve that is formed by the arc from $y_1$ to $y_2$
 in the boundary of $S(x_0)$, plus  the arc from $y_2$ to $y_1$  in the boundary of $S(x_1)$.

Then, we assume the boundary of $S(x_0,...,x_{k-1})$ is a closed curve, when we consider the arc $x_0,...,x_{k-1}, x_{k}$ in $C$, we can prove
the boundary of $S(x_0,...,x_{k})$ is also a closed curve.

We know that we have two closed curves: Suppose that $Q$ is the boundary of $S(x_0,...,x_{k-1})$ , and $R$ is the boundary of $S(x_k)$. $(x_{k-1}, x_k)$ is in
$S(x_k)$, and  $(x_{k-1}, x_k)$ is in $S(x_0,...,x_{k-1})$ . There are two 2-cells $A$, $B$ containing  $(x_k-1, x_k)$ in $S(x_k)\cap S(x_0,...,x_{k-1})$.

$x_{k-1}$ is on the boundary cycle of $S(x_k)$, then $x_{k-1}$ must have two adjacent points in $R$, $y_1$, and $y_2$. $(x_{k-1},y_1)$ and $(x_{k-1},y_2)$ are two edges in $S(x_k)\cap S(x_0,...,x_{k-1})$. In the same way above, we will have the cycle passing $y_1$ and $y_2$ that is the boundary curve of $S(x_0,...,x_{k})$ .

In the rest of the proof, we just need to treat $S(x_0,...,x_m)$ to be $S(x)$ in Case 1. See Fig. 3.5.

\begin{figure}[h]
	\begin{center}

   \epsfbox{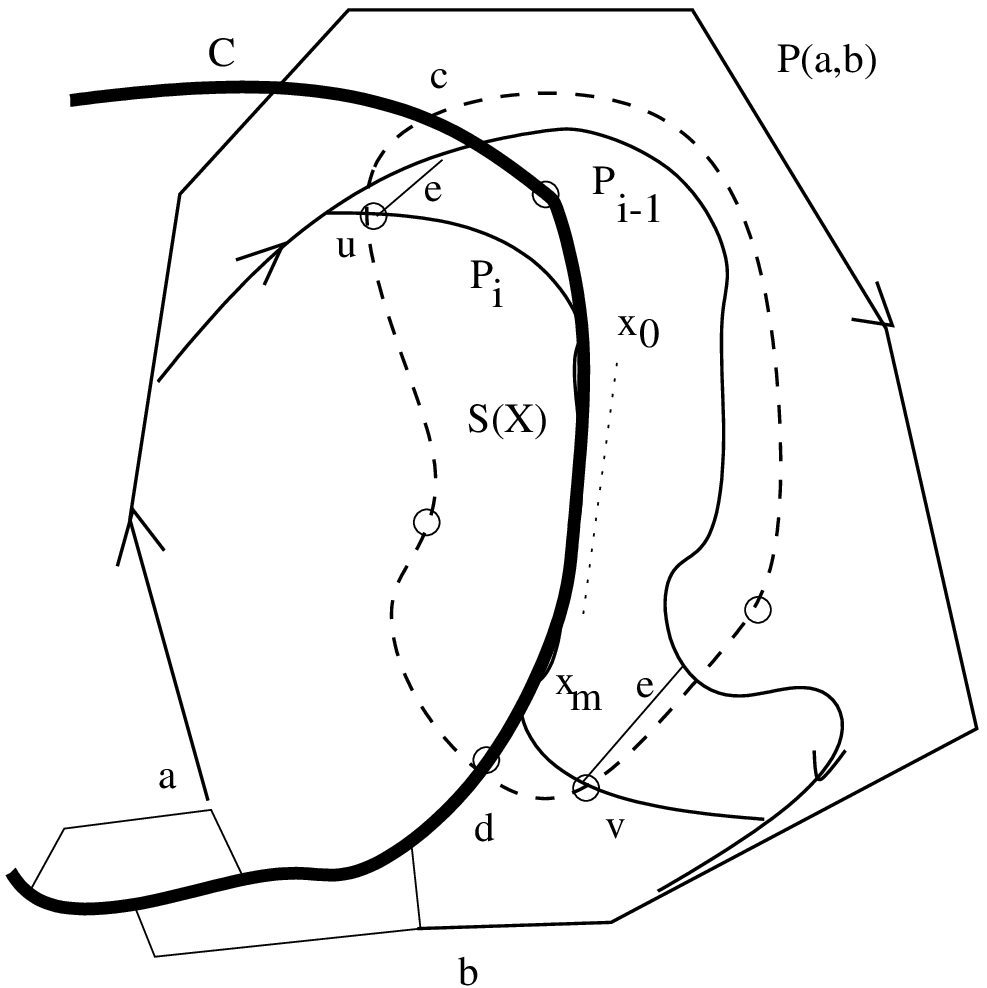 }

\caption*{Figure 3.5  A pair of crossing-over curves pass an arc $X=\{x_0,...,x_m\}$ }
\end{center}

\end{figure}


So we will have,
if $P_{i-1}(a,b)$ and $C$ do not cross over each other, then, $P_{i}(a,b)$ and $C$ will not cross over each other.
Therefore, any $P(a,b)$ must cross over $C$. This completes the proof of the discrete Jordan curve theorem.

Let's first state again that $P_{i}(a,b)$ passes $x_0...x_m$ but $P_{i-1}(a,b)$ does not contain any point of $\{x_0,...,x_m\}$.
In addition, $P_{i-1}(a,b)$  and $P_{i}(a,b)$ is gradually varied, i.e. $P_{i}(a,b)$ was deformed from $P_{i-1}(a,b)$. We also know that
$S(X)=S(x_0,...,x_m)$ is the neighborhood of the arc in $C$ meaning that $x_0,...,x_m$ is a part of the closed curve $C$. The boundary of
$S(X)=S(x_0,...,x_m)$ is a closed curve too.

$u, v,  c,  d$ are on the boundary of  $S(x_0,...,x_m)$  (Assume $u, v,  c,  d$ are not pseudo points, otherwise, we can find
corresponding none-pseudo on  the boundary of  $S(x_0,...,x_m)$.)  $u, (x_0,...,x_m),v$  is a part of $P_n$
We also know that $c$ and $(x_0,...,x_m)$ are not in $P_{n-1}$. There will be two 2-cells,
$U$ and $V$, are in between $P_{i}(a,b)$ and $P_{i-1}(a,b)$ (all points of $U$ and $V$ are in $P_{i}(a,b)\cup P_{i-1}(a,b)$)  such that $(u, x_0) \in U$ and $(x_m,v)\in V$.

Let $P_{n-1}\cap U = S1$ and $P_{n}\cap U = S2$. Let $e1$ be the edge in $U$ linking $S1$ to $S2$ (in most cases, $e1$ incident to $u$, but not necessarily ), and let $e2$ be the edge in $U$ linking $S2$ to $S1$ (possibly starting at $x_0$).
So, $(e2) (S1) (e1) (S2)$ are the boundary of $U$, counterclockwise.

Subcase (i): If $U$ contains $v$ ($U=V$), all points in $U$'s boundary are contained in $S(\{x_0,...,x_m\})$ by the definition of $S(x_0)$.  we will have a cycle $u\cdot d \cdot v (e2) (S1) (e1) $ in the boundary of $S(X=\{x_0,...,x_m\})$
$c$ is on the boundary of $S(X)$ too. But $c\notin P_{i-1}$ It must be in the boundary curves (of $S(X)$)
from $u$ to $d$ or the curve from $d$ to $v$.
Then $c,d$ in $S(X)$ must in the same side of $u X v$ which is part of $P_{i}$.
Therefore, $C$ and $P_{i}(a,b)$ do not cross-over each other at $X$. (See Fig. 3.5.)

Subcase (ii): If $U$ does not contain $v$, then there must be a 2-cell $V$ (in between $P_{i-1}$ and $P_{i}$)
containing $(x_m,v)$.

Let $e1=(p1,p2)$ be the edge in $U$ incident to a point in $P_{i-1}$ and a point in $P_i$, respectively.
(In most cases, $e1$ incident to $u$, i.e. $u=p2$, but not necessarily ). And let $e2=(r2,r1)$ be the edge in $V$ incident to a point in $P_{i}$ and a point in $P_{i-1}$, respectively. $r2$ is usually $v$.

$c$ must not be in $U$, deformation means each point in each 2-cell in between $P_i$ and $P_{i-1}$  must be in $P_i \cup P_{i-1}$. Formally, ($P_i$ $XoRSum$  $P_{i-1}$) is a set of 2-cells; every point in these 2-cells is in $P_n \cup P_{n-1}$ .

We can see that $U$ and $V$ are line-connected in $S(X)$ by the definition of line-connected paths  meaning there is a path of 2-cells where each adjacent pair shares a 1-cell. (See Fig. 3.5 )

From $r1$ to $p1$, there is an arc in $P_{i-1}$ . To prove that all the point in this arc are in the boundary of $S(X)$ we need
to prove each point on the arc must be in a 2-cell that contains a point in $\{ x_0,..., x_m\}$, and this 2-cell is
other than (except this 2-cell is) $U$ or $V$.
It gives us some difficult to prove it.

We found a more elegant way to prove  this case by finding another curve that cross-over $C$.  The method is the following:
If $U\neq V$,  there must be a $x_k$ in $\{ x_0,..., x_m\}$, $x_k$ has an edge linking to $P_{i-1}$. (Otherwise, $u,  x_0,..., x_m, v$
are in a 2-cell that contains some points in  $P_{i-1}$. Therefore, $U=V$.) We can also assume that $k$ is not $m$, otherwise,
$v$ is in $P_{i-1}$, so $U=V$.  See  (See Fig. 3.6 )

\begin{figure}[h]
	\begin{center}

  \epsfxsize=4. in
   \epsfbox{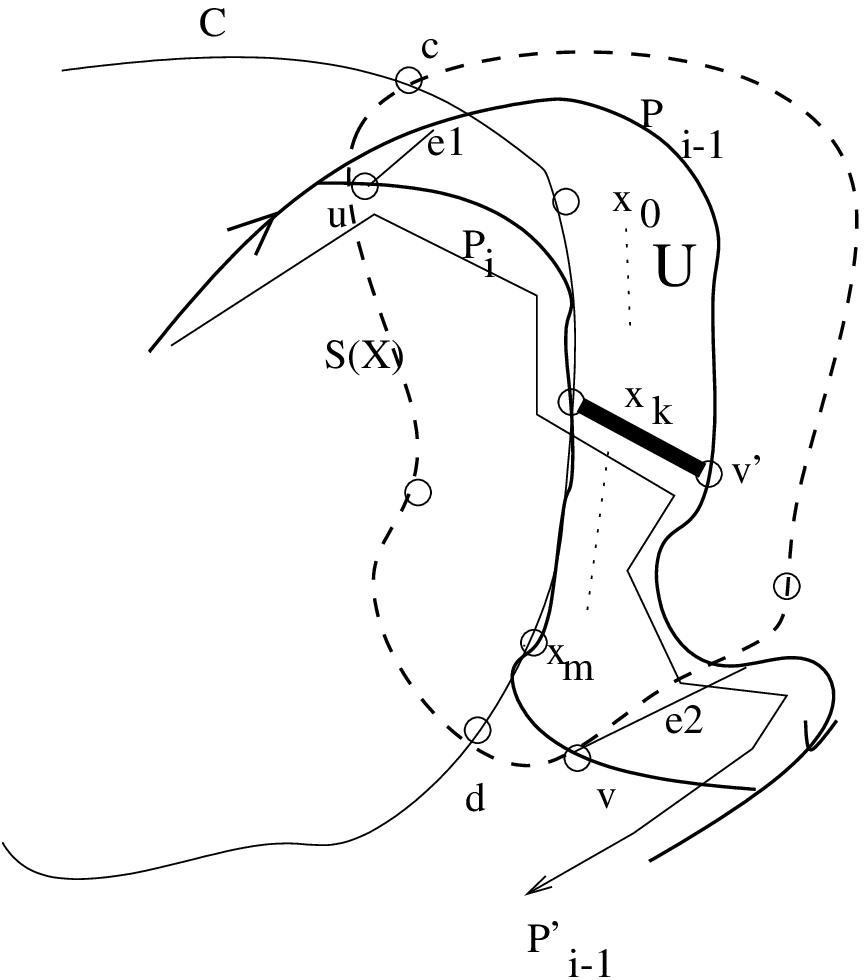 }


\end{center}
\caption*{Figure 3.6  An edge only starts at $x_k$ to $P_{i-1}$;  $x_0,...,x_{k-1}$ do not have any edge to $P_{i-1}$. }

\end{figure}

We select the smallest $k$ having an edge linking to $P_{i-1}$, $0\le k\le m-1$.
Might as well let $(x_k,v')$ is such an edge, and $v's$ is a point in $P_{i-1}$. Therefore, the new curve, $P'_{i}$ that is the same as   $P_{i}$
before and including the point $x_k$, and followed by the partial curve of $P_{i-1}$ after point $v'$ . This curve $P'_{i}=...,u, x_0, ...,x_k, v',...$ does cross-over $C=...,c, x_0, ...,x_k,x_{k+1},...,x_{m}, d,... $ .  It is obvious that $P_{i-1}, P'_{i}, P_{i}$ are gradually varied.
We just inserted a path in between of $P_{i-1}$ and $P_{i}$. This new path $P'_{i-1}$ has such a good property that is
$\{x_0, ...,x_{k-1}\}$ do not have an edge in $P_{i-1}$ . Since $v'$ is in $P'_{i-1}$,  the 2-cell $V'$ (in $P_{i-1}$ $XoRSum$ $P'_{i-1}$) contains $v'$ also contains $(x_{k-1} x_k)$ and $(x_k, v)$ (in $S(x_k)$). Since no edge from $x_0,...,x_{k-1}$ to $P_{i-1}$ , $U$ containing $u$ is just $V'$.
We will have just Subcase (i).

The entire theorem is proven.   $\diamondsuit$\\

Theorem 3.1, the discrete Jordan curve theorem, has a little difference from the classical description of The Jordan curve Theorem.
This is because that discrete curve has its own strict property: $C$ does not contain any 2-cell. In order to satisfy the classical
form. We need to use central pseudo points for each type of cells, especially 1-cells (line-cells) and 2-cells (surface-cells)
So we will allow the simple path (semi-curve) in the prove of  Theorem 3.1. (A little modification will do the task in the following
theorem.) The rest of work is to  prove that there are only two (connected) components in in $S-C$. \\

{\bf Theorem 3.2} (The Jordan Curve Theorem for Generalized Curves) Let $S$ be a discrete simply connected surfaces, ($S$ can be closed or a discrete plane embedded
in 2D Euclidean Space). A closed simple path (0-cell connected semi-curve) $C$ which does not contain any point of $\partial S$
divides $S$ into two components (in terms of allowing central pseudo points for each cell). In other words, $S-C$ consists of two components. These two components are disconnected.

{\bf Proof}  We have now the central pseudo points for each  1-cells and 2-cells. In the proof of Theorem 3.1. We know that
$2-cell$ $a\in A$ and $b\in B$ are not connected in $S-C$. $C$ has orientation of counterclockwise. $(p,r)$ is  counterclockwise
in $A$, and $(p,r)$ is clockwise in $B$. We call $A$ is counterclockwise, and $B$ is  clockwise.  All 2-cells that has an edge $e_i=(p,r)$ in $C$ will have two 2-cells  containing $e_i$, $A_i$
and $B_i$. We always assume that $A_i$ is counterclockwise and  $B_i$ is clockwise. We now add all the central pseudo points to $S$ and remove
the central pseudo points from $C$. Since each 2-cell must have at least three boundary 1-cells.

{\bf Case 1:} A simple path could be just the boundary of
a 2-cell. In this case, we have a central point in the cell, denoted $A$. In addition, the rest of $S-C$ is point-connected. This is because that
each other cell has an edge not in $C$ shared by two 2-cells. Those points in two cells are connected including pseudo points. ( a central pseudo
point in 2-cell always connected to the points on its boundary points.) There are finite number of 2-cells. This process will stop. If there
are two components in $S-A$, a component must have a boundary, this boundary must contain at least an edge $e$ that is not in $\partial S$ and $e$ is not in $C$. Because $C$ does not include any point in $\partial S$. $e$ will be contained by two cells, the center pseudo points of the two
cells are connected by the center pseudo point of $e$. $e$ is not a boundary edge of the set. Therefore the boundary is only $C$ or $\partial S$

{\bf Case 2:} Let's prove the case of there are more 2-cells in counterclockwise, i.e. $A_i$, $i\ge 1$. We can prove that all $A_i$ are connected. This is because that any point $p$ in $C$ is contained by two 1-cells $e1$ and $e2$ in $C$. These two 1-cells are contained by $A_i$ and $A_j$,
respectively. If $A_i$ and $A_j$ share an edge, then, the central pseudo points of $A_i$ and $A_j$ are connected.
If  $A_i$ and $A_j$ do not share an edge, we know $A_i$ and $A_j$ are in $S(p)$, there must be a cycle contains some edges in $A_i$ and
some edges of $A_j$, and $e1 \cup e2$. So $A_i$ and $A_j$ are connected (meaning their central pseudo points ) do not pass $e1 \cup e2$.
Therefore, all $A_i$'s (meaning their central pseudo points) are connected.

(All 2-cells sharing an edge $e\in A_{i}$, $e$ is not in $C$, are
connected. Since $S$ is simply connected, there will be a (side-)gradually varied sequence of closed paths from $C$ to boundary of $A_i$.
All cells pass through the sequence are connected. There is only one component. Seems not necessary since we have Theorem 3.1.)

In the same way, we can prove that all $B_i$'s are connected. Extend $B_i$ by connecting its edge not in $C$ will result several components.
Those components does not have an edge as we discussed above (in case 1) except $C$ and $\partial S$. So the component that contain
a $B_i$ will be just a component.

Let's now prove that any point $p$ in $S-C$, must be contained in the component containing $A_i$ or the component containing $B_i$.
We know that any two points are point-connected by a path in $S$. Let $c\in C$, $P(p,c)$ is such a path. There must be a first point
in $P(p,c)$, $p'$, that is adjacent a point $c'\in C$ ( may or may not be point $c$). $(p', c')$ must belong to an $A_i$ or $B_j$. So
If $(p', c')$ belong to $A_i$, it is point connected to the central pseudo points of $A_i$ . We call it component $A$.  All points in $A$
are connected since $A_i$ are connected for all $i$.

If $(p', c')$ belong to $B_j$, it is point connected to the central pseudo points of $B_j$ . We call it component $B$. All points in $B$
are connected since $B_j$ are connected for all $j$.

Points in $A_i-C$ is not connected to $B_j$ in $S-C$ based on Theorem 3.1.  Therefore, any point in $A$ is not connected to any point in $B$ in
in $S-C$. We now complete the proof of Theorem 3.2, the general Jordan Curve Theorem.
$\diamondsuit$\\

So we can allow the simple path (semi-curve). This is the general case of Jordan Curve Theorem.  \\\\

\centerline{\bf Discussion}

We have added more details to the original proof of  the discrete Jordan curve theorem.
The idea of proof is unchanged.

Up to now, we can see that we have proved completely the discrete Jordan curve theorem.
If we embed the discrete surface into a plane. The proof is still valid.

This proof is general and extendable not only for Euclidean space. There is no
approximation process since the boundary of a simple closed surface can be grabbed
as any type of shapes. \\\\

\centerline{\large {References}}
\begin{enumerate}

\item [1] G.T. Herman, {\it Geometry of Digital Spaces,} Birkhauser, Boston, 1998.

\item [2]  E. Artzy, G. Frieder and G. T. Herman, ``The theory, design,
	 implementation and evaluation
	 of a three-dimensional surface detection algorithm,''
	 {\it Comput. Vision Graphics Image Process.}  Vol 15 , pp. 1-24, 1981.

\item [3] D. G. Morgenthaler and A. Rosenfeld, ``Surfaces in three-dimensional images,''
	 {\it Inform. and Control,} Vol 51, pp. 227-247, 1981.

\item [4] L. Chen, ``Generalized discrete object tracking algorithms and
	implementations,'' In Melter, Wu, and Latecki ed, {\it Vision Geometry VI},
	SPIE Vol. 3168, pp 184-195, 1997.

\item [5]  L. Chen, ``Generalized Discrete Object (I): Curves, Surfaces and Manifolds,''
	 Manuscript,  1998.

\item[6]  L. Chen, ``Point spaces and raster spaces in digital geometry and topology,''
	 in Melter, Wu, and Latecki ed, {\it Vision Geometry VII},SPIE Proc.
	 3454, pp 145-155,1998.

\item [7]  L. Chen, ``(alpha, beta)-type digital surfaces and general digital
	surfaces,'' in Melter, Wu, and Latecki ed, {\it Vision Geometry VII},
  SPIE Proc. 3454, pp 28-39,1998.

\item [8]   L. Chen, H. Cooley and J. Zhang, ``The equivalence between two
  definitions of digital surfaces,'' {\it Information Sciences}, Vol 115, pp
  201-220, 1999.

\item [9] T.Y.Kong and  A.Rosenfeld, ``Digital topology: Introduction and survey,''
	  {\it Comput. Vision Graphics Image Process.} Vol 48, pp. 357-393, 1989.

\item [10]  L. Chen and J. Zhang, ``Digital manifolds: A Intuitive Definition and Some
	 Properties'', {\it The Proc. of the Second ACM/SIGGRAPH Symposium on Solid
	 Modeling and Applications,} pp. 459-460,  Montreal, 1993.

\item [11]  L. Chen and J. Zhang, ``Classification of Simple digital Surface Points
	 and A Global Theorem for Simple Closed Surfaces'',  in Melter and Wu ed,
	 {\it Vision Geometry II}, SPIE Vol 2060, pp. 179-188, 1993.

\item [12] S. Lefschetz, ``Introduction to Topology,'' Princeton University Press
			  New Jersey, 1949.

\item [13]  A. Rosenfeld, ```Continuous' functions on digital pictures,'' {\it Pattern
	  Recognition Letters}, No 4, pp 177-184, 1986.

\item [14] Z. Pawlak, ``Rough calculus,'' in P.P Wang ed, {\it Advances in Machine
	Intelligence $\&$ Soft-Computing}, Vol IV,
		 Duke University, 1996.

\item [15] L. Chen, ``Three-dimensional fuzzy digital topology and its
	  applications(I),''
		{\it Geophysical Prospecting for petroleum}, Vol 24, pp 86-89, 1985.

\item [16] L. Chen, ``The necessary  and  sufficient  condition  and  the
		efficient algorithms for gradually varied fill,'' {\it Chinese Science
		Bulletin}, Vol 35,10(1990). (or {\it Abstracts of SIAM Conference on
		Geometric design,} Temple, AZ, 1989.)

\item [17]  L. Chen, ``Gradually varied surface and its optimal uniform approximation'',
	 {\it IS\&T/SPIE Symposium on Electronic Imaging,} SPIE Proc.
	 Vol 2182, pp. 300-307, 1994.

\item [18]  M. Newman, {\it Elements of the Topology of Plane Sets of Points,}
	  Cambridge, London, 1954.

\item [19]  L. Latecki, "3D well-composed pictures",  {\it The Proc.
	 of SPIE on Vision Geometry IV,} Vol 2573,  pp. 196-203, 1995.

\item [20]  L. Latecki, U. Eckhardt, and A. Rosenfeld, "3D well-composedness of
	  digital sets",   {\it The Proc.  of SPIE on Vision Geometry II,}
		Vol 2060,  pp. 61-68, 1993.

\item [21] T. H. Cormen, C.E. Leiserson, and R. L. Rivest,  {\it Introduction to
	 Algorithms,} MIT Press, 1993.

\item [22] O. Veblen, Theory on Plane Curves in Non-Metrical Analysis Situs, {\it Transactions of the American Mathematical Society } 6 (1): 83–98, 1905.

\item [23]  W.T. Tutte,  Combinatorial oriented maps. Can. J. Math. XXXI:5 (1979), 986–1004.

\item [24] J.-F. Dufourd, An intuitionistic proof of a discrete form of the Jordan curve theorem
           formalized in Coq with combinatorial hypermaps, Journal of Automated Reasoning 43 (1) (2009) 19–51.

\item [25] L. Chen, Discrete Surfaces and Manifolds, SP Computing, Rockville, 2004.

\end{enumerate}

\end{document}